\documentclass[12pt]{amsart}
\usepackage{amscd,amssymb,amsmath,multicol}
\newtheorem{thm}[equation]{Theorem}
\numberwithin{equation}{section}
\newtheorem{cor}[equation]{Corollary}

\newtheorem{lem}[equation]{Lemma}

\newtheorem{defin}[equation]{Definition}

\newtheorem{prop}[equation]{Proposition}
\newtheorem{tab}[equation]{Table}

\begin{document}
\raggedbottom \voffset=-.7truein \hoffset=0truein \vsize=8truein
\hsize=6truein \textheight=8truein \textwidth=6truein
\baselineskip=18truept

\def\ss{\smallskip}
\def\ssum{\sum\limits}
\def\dsum{{\displaystyle{\sum}}}
\def\la{\langle}
\def\ra{\rangle}
\def\on{\operatorname}
\def\a{\alpha}
\def\bz{{\Bbb Z}}
\def\eps{\epsilon}
\def\br{{\bold R}}
\def\bc{{\bold C}}
\def\bN{{\bold N}}
\def\nut{\widetilde{\nu}}
\def\tfrac{\textstyle\frac}
\def\b{\beta}
\def\G{\Gamma}
\def\g{\gamma}
\def\zt{{\bold Z}_2}
\def\bs{{\bold s}}
\def\bg{{\bold g}}
\def\bof{{\bold f}}
\def\bq{{\bold Q}}
\def\be{{\bold e}}
\def\line{\rule{.6in}{.6pt}}
\def\xb{{\overline x}}
\def\xbar{{\overline x}}
\def\ybar{{\overline y}}
\def\zbar{{\overline z}}
\def\ebar{{\overline \be}}
\def\nbar{{\overline n}}
\def\rbar{{\overline r}}
\def\Ubar{{\overline U}}
\def\et{{\widetilde e}}
\def\ni{\noindent}
\def\ms{\medskip}
\def\ehat{{\hat e}}
\def\nbar{{\overline{n}}}
\def\minp{\min\nolimits'}
\def\N{{\Bbb N}}
\def\Z{{\Bbb Z}}
\def\Q{{\Bbb Q}}
\def\R{{\Bbb R}}
\def\C{{\Bbb C}}
\def\mo{\on{mod}}
\def\dstyle{\displaystyle}
\def\Remark{\noindent{\it  Remark}}
\title
{Divisibility by 2 of partial Stirling numbers}
\author{Donald M. Davis}
\address{Department of Mathematics, Lehigh University\\Bethlehem, PA 18015, USA}
\email{dmd1@lehigh.edu}
\date{September 22, 2011}


\maketitle
\begin{abstract} The partial Stirling numbers $T_n(k)$ used here are defined as $\dstyle{\sum_{i\text{ odd}}\tbinom nii^k}$. Their
2-exponents $\nu(T_n(k))$ are important in algebraic topology. We provide many specific results, applying to all values of $n$,
stating that, for all $k$ in a certain congruence class mod $2^t$, $\nu(T_n(k))=\nu(k-k_0)+c_0$, where $k_0$ is a 2-adic integer and $c_0$ a positive
integer. Our analysis involves several new general results for $\nu(\sum\binom n{2i+1}i^j)$, the proofs of which involve a new family of
polynomials. Following Clarke (\cite{Cl}), we interpret $T_n$ as a function on the 2-adic integers, and the 2-adic integers $k_0$ described
above as the zeros of these functions.
\end{abstract}
\section{Main results}\label{intro}
The partial Stirling numbers $T_n(k)$ used here are defined, for integers $n$ and $k$ with $n$ positive, by
$$T_n(k)=\sum_{i\text{ odd}}\tbinom ni i^k.$$
Other versions can be defined localized at other primes and summed over restricted congruences.
Let $\nu(-)$ denote the exponent of 2 in an integer. The numbers $\nu(T_n(k))$ are important in
algebraic topology (\cite{BDSU}, \cite{CK}, \cite{D23}, \cite{DP}, \cite{DS1}, \cite{Lune}), and
work on evaluating these numbers has appeared in the above papers as well as \cite{Cl}, \cite{D2},  \cite{Lun},  \cite{DS2}, and \cite{You}.
In this paper, we give complete results for $n\le 36$ and also for $n=2^e+1$ and $2^e+2$, and we give two families of results applying to all values of $n$ but with $k$ restricted to
certain congruence classes. In \cite{SU2e}, some of these results will be applied to obtain new results for $v_1$-periodic homotopy
groups of the special unitary groups. We also present in Section \ref{thmpfsec} some new results about $\nu(\sum_i\tbinom n{2i+1}i^k)$.
The proofs of these, in Section \ref{sumpfsec}, introduce a new family of polynomials $q_m(x)$, which might be of independent interest. Finally, in Section \ref{Clsec} we discuss our
results in the context of analytic functions on the 2-adic integers, and Hensel's Lemma.

We begin with the result which is easiest to state, and hence best illustrates the nature of our results.
\begin{thm}\label{thm0} Let $e\ge2$, $n=2^e+1$ or $2^e+2$, and $1\le i\le 2^{e-1}$.
\begin{enumerate}
\item There is a 2-adic integer $x_{i,n}$ such that for all integers $x$
$$\nu(T_n(2^{e-1}x+i))=\nu(x-x_{i,n})+n-2.$$
Moreover
$$x_{i,{2^e+1}}\equiv\begin{cases}1+2^i\pmod{2^{i+1}}&\text{if $i=2^{e-2}$ or $2^{e-1}$}\\
1\pmod{2^{i+1}}&\text{otherwise}\end{cases}$$
and
$$x_{i,{2^e+2}}\equiv\begin{cases}1+2^{i-1}\pmod{2^i}&\text{if }1\le i\le 2^{e-2}\\
1+2^i\pmod{2^{i+1}}&\text{if }2^{e-2}<i<2^{e-1}\\
1\pmod{2^{i+1}}&\text{if }i=2^{e-1}.\end{cases}$$
\item Let $g(x)=\nu(T_n(2^{e-1}x+i))-(n-2)$.
Then $x_{i,n}=2^{t_0}+2^{t_1}+\cdots$, where $t_0=g(0)$ and $t_{j+1}=g(2^{t_0}+\cdots+2^{t_j})$.
\end{enumerate}\end{thm}
\noindent If $e=1$, the result is different.  See Table \ref{bigtbl}.

  For example, the last 24 digits of the binary expansion of $x_{4,9}$ are $$100000000001101010110001,$$
and so we can make the following more explicit statement.
$$\nu(T_9(4x+4))=\begin{cases}7&x\equiv0\ (2)\\
8&x\equiv3\ (4)\\
9&x\equiv5\ (8)\\
10&x\equiv9\ (16),\end{cases}$$
and continue indefinitely, just noting the last position in which the binary expansions of $x$ and $x_{4,9}$ differ.

Our next result utilizes {\tt Maple} calculations in its proof. Although each case applies to infinitely many values of $x$, we will explain in the proof how
each case can be reduced to a small number of verifications.
\begin{thm}\label{36} For each $n\le 36$, there is a partition of $\Z$ into finitely many congruence classes $C=[i\mod 2^m]$ such that, for each, either (a) there exists
a 2-adic integer $x_0$  and a positive integer $c_0$ such that $\nu(T_n(2^mx+i))=\nu(x-x_0)+c_0$ for all integers $x$, or (b)
there exists a positive integer $y_0$ such that $\nu(T_n(k))=y_0$ for all $k$ in $C$. The congruence classes $C$ and integers $c_0$ and $y_0$ are as in Tables \ref{bigtbl} and \ref{bigtbl2}.

Let $g(x)=\nu(T_n(2^mx+i))-c_0$. Then $x_0=2^{t_0}+2^{t_1}+\cdots$, with $t_0=g(0)$
and $t_{j+1}=g(2^{t_0}+\cdots+2^{t_j})$.
\end{thm}

 We conjecture that the general form of the theorem
can be extended to all integers $n$; i.e., that for each $n$ there is a partition of $\Z$ into finitely many congruence classes on each of which
either $\nu(T_n(k))=\nu(k-k_0)+c_0$ for some $k_0$ and $c_0$ or else $\nu(T_n(k))$ is constant on $C$. In the tables, the letter $i$ refers
to any integer.

\medskip
\medskip
\begin{minipage}{6.5in}
\begin{tab}\label{bigtbl}{Values of $C$, $c_0$, and $y_0$ in Theorem \ref{36}}
\begin{center}
\begin{tabular}{c|c|c|c|c|c|c|c|c}
$n$&$C$&$c_0$&$y_0$&\qquad\qquad\qquad&$n$&$C$&$c_0$&$y_0$\\
\hline
$3$&$0\ (2)$&&$2$&&$4$&$0\ (2)$&&$3$\\
&$1\ (2)$&&$1$&&&$1\ (2)$&&$4$\\
\hline
$5$&$0,1\ (2)$&$3$&&&$6$&$0,1\ (2)$&$4$&\\
\hline
$7$&$0\ (2)$&$6$&&&$8$&$0\ (2)$&$7$&\\
&$1\ (2)$&$4$&&&&$1\ (2)$&$9$&\\
\hline
$9$&$i\ (4)$&$7$&&&$10$&$i\ (4)$&$8$&\\
\hline
$11$&$0,2\ (4)$&$9$&&&$12$&$0,2\ (4)$&$10$&\\
&$1,3\ (4)$&$8$&&&&$1,3\ (4)$&$11$&\\
\hline
$13$&$1,2\ (4)$&$10$&&&$14$&$2,3\ (4)$&$11$&\\
&$0,4\ (8)$&$12$&&&&$0,1\ (8)$&$13$&\\
&$7\ (8)$&&$11$&&&$4,5\ (8)$&$13$&\\
&$3\ (16)$&&$13$&&&&&\\
&$11\ (16)$&&$15$&&&&&\\
\hline
$15$&$3\ (4)$&$11$&&&$16$&$0\ (4)$&$15$&\\
&$0\ (4)$&$14$&&&&$1\ (4)$&$18$&\\
&$1,5\ (8)$&$13$&&&&$2,6\ (8)$&$17$&\\
&$2,6\ (8)$&$16$&&&&$3,7\ (8)$&$20$&\\
\hline
$17$&$i\ (8)$&$15$&&&$18$&$i\ (8)$&$16$&\\
\hline
$19$&$0,2,4,6\ (8)$&$17$&&&$20$&$0,2,4,6\ (8)$&$18$&\\
&$1,3,5,7\ (8)$&$16$&&&&$1,3,5,7\ (8)$&$19$&\\
\hline
$21$&$1,2,5,6\ (8)$&$18$&&&$22$&$0,1\ (8)$&$20$&\\
&$0,3\ (8)$&$19$&&&&$2,3,6,7\ (8)$&$19$&\\
&$7,15\ (16)$&$21$&&&&$4,12\ (16)$&$22$&\\
&$12\ (16)$&&$20$&&&$5,21\ (32)$&$24$&\\
&$4\ (32)$&&$22$&&&$13\ (16)$&&$21$\\
&$52\ (64)$&&$24$&&&&&\\
&$84\ (128)$&&$26$&&&&&\\
&$20\ (256)$&&$28$&&&&&\\
&$148\ (256)$&&$29$&&&&&\\
\hline
\end{tabular}
\end{center}
\end{tab}
\end{minipage}
\medskip
\newpage

\medskip
\begin{minipage}{6.5in}
\begin{tab}\label{bigtbl2}{More values of $C$, $c_0$, and $y_0$ in Theorem \ref{36}}
\begin{center}
\begin{tabular}{c|c|c|c|c|c|c|c|c}
$n$&$C$&$c_0$&$y_0$&\qquad\qquad\qquad&$n$&$C$&$c_0$&$y_0$\\
\hline
$23$&$0,4\ (8)$&$21$&&&$24$&$0,4\ (8)$&$22$&\\
&$3,7\ (8)$&$19$&&&&$1,5\ (8)$&$24$&\\
&$1\ (8)$&$20$&&&&$2\ (8)$&$23$&\\
&$2\ (8)$&$22$&&&&$3\ (8)$&$25$&\\
&$6,22\ (32)$&$26$&&&&$6,22\ (32)$&$27$&\\
&$5,21\ (32)$&$24$&&&&$7,23\ (32)$&$29$&\\
&$13\ (16)$&&$21$&&&$14\ (16)$&&$24$\\
&$14\ (16)$&&$23$&&&$15\ (16)$&&$26$\\
\hline
$25$&$1,2,3,4\ (8)$&$22$&&&$26$&$2,3,4,5\ (8)$&$23$&\\
&$5,6,7,8\ (16)$&$24$&&&&$0,1,7,8,9,15\ (16)$&$25$&\\
&$0,13,14,15\ (16)$&$24$&&&&$6,22\ (32)$&$27$&\\
&&&&&&$14\ (16)$&&$24$\\
\hline
$27$&$3,5\ (8)$&$23$&&&$28$&$4,6\ (8)$&$25$&\\
&$4,6\ (8)$&$24$&&&&$5,7\ (8)$&$26$&\\
&$1,7,9,15\ (16)$&$25$&&&&$0,2,8,10\ (16)$&$27$&\\
&$0,2,8,10\ (16)$&$26$&&&&$1,3,9,11\ (16)$&$28$\\
\hline
$29$&$5,6\ (8)$&$25$&&&$30$&$6,7\ (8)$&$26$&\\
&$7\ (8)$&$26$&&&&$2,3,10,11\ (16)$&$28$&\\
&$1,2,9,10\ (16)$&$27$&&&&$0,1,8,9\ (16)$&$29$&\\
&$0,4,8,12\ (16)$&$28$&&&&$4,5,12,13\ (16)$&$29$&\\
&$11\ (16)$&$29$&&&&&&\\
&$3\ (32)$&&$30$&&&&&\\
&$19\ (32)$&&$31$&&&&&\\
\hline
$31$&$7\ (8)$&$26$&&&$32$&$0\ (8)$&$31$&\\
&$0\ (8)$&$30$&&&&$1\ (8)$&$35$&\\
&$3,11\ (16)$&$28$&&&&$4,12\ (16)$&$33$&\\
&$1,5,9,13\ (16)$&$29$&&&&$2,6,10,14\ (16)$&$34$&\\
&$4,12\ (16)$&$32$&&&&$5,13\ (16)$&$37$&\\
&$2,6,10,14\ (16)$&$33$&&&&$3,7,11,15\ (16)$&$38$&\\
\hline
$33$&$i\ (16)$&$31$&&&$34$&$i\ (16)$&$32$&\\
\hline
$35$&$2i\ (16)$&$33$&&&$36$&$2i\ (16)$&$34$&\\
&$2i+1\ (16)$&$32$&&&&$2i+1\ (16)$&$35$&
\end{tabular}
\end{center}
\end{tab}
\end{minipage}
\medskip

One can notice a lot of nice patterns in these tables, and formulate (and sometimes prove)
conjectures about their extension to all values of $n$. One interesting idea, following Clarke (\cite{Cl}),
is to note that since $T_n(k)$ mod $2^m$ only depends on $k$ mod $2^{m-1}$,  $T_n(-)$ extends to a function
$T_n:\zt\to\zt$, where $\zt$ denotes the 2-adic integers. Here the metric on $\zt$ is given, as usual, by $d(x,y)=|x-y|$,
where $|z|:=1/2^{\nu(z)}$. The 2-adic integers $2^mx_0+i$ which occur in Theorem \ref{36}
 are just the zeros of the function $T_n$. We can count the number of zeros to be given
as in Table \ref{ztbl}, and might try to formulate a guess about the general formula for this number of zeros.

\medskip
\begin{minipage}{6.5in}
\begin{tab}\label{ztbl}{Number of zeros of $T_n$}
\begin{center}
\begin{tabular}{c|c}
$n$&Number of 0's of $T_n$\\
\hline
$1$-$4$&$0$\\
$5$-$8$&$2$\\
$9$-$13$&$4$\\
$14$-$16$&$6$\\
$17$-$21$&$8$\\
$22$-$24$&$10$\\
$25$-$29$&$12$\\
$30$-$32$&$14$\\
$33$-$36$&$16$
\end{tabular}
\end{center}
\end{tab}
\end{minipage}
\medskip

Our second general result establishes for all $n$, except those 1 less than a 2-power, the values of $\nu(T_n(k))$
for $k$ in the congruence class containing $0$. We could almost certainly include $n=2^{e}-1$ into this theorem,
but the details of proving that case are so detailed as to be perhaps not worthwhile here.
\begin{thm} Let $n\ge5$ and
$$(a,b)=\begin{cases}(-2,1)&\text{if }n=2^{e}\\
(-1,2)&\text{if }2^{e}<n\le 3\cdot2^{e-1}\\
(0,1)&\text{if }3\cdot2^{e-1}<n\le2^{e+1}-2.
\end{cases}$$
Then there exists a $2$-adic integer $\xb_n$ such that for all integers $x$
$$\nu(T_n(2^{e+a}x)) =\nu(x-\xb_n)+n-b.$$
\label{thm2}
\end{thm}
\noindent The cases $n=2^e+1$ and $2^e+2$ of this theorem overlap with Theorem \ref{thm0}.
For these $n$, we have $\xb_{n}=1+x_{2^{e-1},n}$. For all $n$ in Theorem \ref{thm2}, there is an algorithm
for $\xb_n$ totally analogous to that of Theorem \ref{thm0}.

Our next result is of a similar nature, but applies to many more congruence classes.
The cases to which it applies are those in which the 2-exponent of a certain sum (see (\ref{t3e}) and (\ref{t3e2})) is determined
by exactly one of its summands, and for which the mod 4 result \ref{refine} suffices to prove it. The  algorithm for computing $x_0$ is like that of Theorem \ref{36}.
Here and throughout, $\a(n)$ denotes the number of 1's in the binary expansion of $n$.

\begin{thm}\label{thm3} Suppose $2^e+2^t\le n<2^e+2^{t+1}$ with $0\le t\le e-1$. Let
$$S_n=\{p:\ \max(0,n-2^e-2^{e-1})\le p<2^{e-1}\text{ and }\tbinom{n-1-p}p\equiv1\ (2)\}.$$
If $p\in S_n$, say that an integer $q<2^{e-1}$ is {\em associated} to $p$ if $q=p$ or
$q=p+2^{w}$
with $w=\nu(n)-1$ or $w>t$. If $q$ is associated to an integer $p$ of $S_n$, then there
exists a $2$-adic integer $x_0$ such that for all integers $x$
$$\nu(T_n(2^{e-1}x+q))=\nu(x-x_0)+n-2-\a(p_0),$$
where $p_0$ is the residue of $p$ mod $2^t$.
\end{thm}

  A bit of work is required to get any sort of feel for the complicated condition in this theorem.
In Table \ref{thm3tbl}, we list for $n$ from 17 to 31, the values of $p$ in $S_n$, then the additional values of $q$ covered by the theorem,
 and finally the values of $i$ for which Theorem \ref{36}, as depicted in Tables \ref{bigtbl} and \ref{bigtbl2}, gives a value for the
 congruence $i$ mod 8 which is not covered by Theorem \ref{thm3}. The strength of the theorem is, of course, that it applies to all values of $n$ (except 2-powers).

\medskip
\begin{minipage}{6.5in}
\begin{tab}\label{thm3tbl}{Comparison of Theorems \ref{thm3} and \ref{36}}
\begin{center}
\begin{tabular}{c|lll}
$n$&$p\in S_n$&additional $q$&mod $8$ results missed\\
\hline
$17$&$0,1,2,4$&$3,5,6$&$7$\\
$18$&$0,2,4$&$1,3,5,6$&$7$\\
$19$&$0,1,3,4,5$&$7$&$2,6$\\
$20$&$0,4$&$2,6$&$1,3,5,7$\\
$21$&$0,1,2,5,6$&&$3$\\
$22$&$0,2,6$&$1,3,7$&\\
$23$&$0,1,3,7$&&$2,4$\\
$24$&$0$&$4$&$1,2,3,5$\\
$25$&$1,2,4$&&$3$\\
$26$&$2,4$&$3,5$&\\
$27$&$3,4,5$&&$6$\\
$28$&$4$&$6$&$5,7$\\
$29$&$5,6$&&$7$\\
$30$&$6$&$7$&\\
$31$&$7$&&$0$\\
\end{tabular}
\end{center}
\end{tab}
\end{minipage}

\section{Proofs of main theorems}\label{thmpfsec}
In this section, we prove the four main theorems listed in Section \ref{intro}.  A central ingredient in the proofs is results
about $\nu\bigl(\sum_i\binom n{2i+1}i^k\bigr)$. We begin by providing six results about this, of which all but
the first are new. The proofs of most of these appear in Section \ref{sumpfsec}.

The first result was proved in \cite[3.4]{DS1}.
\begin{prop}\label{DSthm} $($\cite[3.4]{DS1}$)$ For any nonnegative integers $n$ and $k$,
$$\nu\bigl(\sum_i\tbinom n{2i+1}i^k\bigr)\ge\nu([n/2]!).$$
\end{prop}
In using this, and many times throughout the paper, we use
\begin{equation}\label{ale}\nu(n!)=n-\a(n).\end{equation}
The next result is a refinement of Proposition \ref{DSthm}. Here and throughout, $S(n,k)$ denote Stirling numbers of the second kind.
\begin{prop}\label{stirowr} Mod $4$
$${\tfrac 1{n!}}\sum_i\tbinom{2n+\eps}{2i+b}i^k\equiv\begin{cases}S(k,n)+2nS(k,n-1)&\eps=0,\,b=0\\
(2n+1) S(k,n)+2(n+1)S(k,n-1)&\eps=1,\,b=0\\
2nS(k,n-1)&\eps=0,\,b=1\\
S(k,n)+2(n+1)S(k,n-1)&\eps=1,\,b=1.
\end{cases}$$
\end{prop}

The proofs of the last three propositions all involve new polynomials $q_m(x)$, which might be of independent interest.
See Definition \ref{qdef} for the definition, which pervades Section \ref{sumpfsec}.
\begin{prop}\label{qprop1} For any nonnegative integers $n$ and $k$,
$$\nu\bigl(\sum_i\tbinom n{2i+1}i^k\bigr)\ge n-k-\a(n).$$
\end{prop}

\begin{prop}\label{qprop2} For any nonnegative integers $n$ and $k$ with $n>k$,
$$\nu\bigl(\sum_i\tbinom n{2i+1}i^k\bigr)\ge n-1-k-\a(k)$$
with equality iff $\binom{n-1-k}k$ is odd.
\end{prop}

The final proposition is a refinement of Proposition \ref{qprop2}.
\begin{prop}\label{refine}  If $n$ and $k$ are nonnegative integers with $n>k$, then, mod 4,
$$\sum_i\tbinom n{2i+1}i^k/(2^{n-1-2k}k!)\equiv\tbinom{n-1-k}k+\begin{cases}2\tbinom{n-1-k}{k-2}&\text{if $n-1$ and $k$ are even}\\
0&\text{otherwise.}\end{cases}$$
\end{prop}

The following corollary will also be useful.
\begin{cor}\label{lcor1} For $n\ge3$, $j>0$, and $p\in\Z$, $$\nu(\sum\tbinom{n}{2i+1}(2i+1)^pi^j)\ge \max(\nu([\tfrac n2]!),n-\a(n)-j)$$ with equality if $n\in\{2^e+1,2^e+2\}$ and $j=2^{e-1}$. \end{cor}
\begin{proof} The sum equals $\sum_{k\ge0}T_k$, where
$$T_k=2^k\tbinom pk\sum_i\tbinom n{2i+1}i^{j+k}.$$
By Proposition \ref{DSthm}, $\nu(T_k)\ge\nu([\frac n2]!)$, while by Proposition \ref{qprop1}, $\nu(T_k)\ge n-\a(n)-j$, implying the desired inequality.
If $n=2^e+1$ and $j=2^{e-1}$, then $\nu(T_0)=2^{e-1}-1$ by Proposition \ref{qprop2}, while for $k>0$, $\nu(T_k)\ge 2^{e-1}$ by \ref{DSthm}.
If $n=2^e+2$ and $j=2^{e-1}$, $\nu(T_0)=2^{e-1}$ by \ref{qprop2}, $\nu(T_1)>2^{e-1}$ by \ref{stirowr}, and $\nu(T_k)>2^{e-1}$ for $k>1$ by \ref{DSthm}.
\end{proof}

Our proofs of the theorems of Section \ref{intro} will make essential use of the following result of
\cite{D2}. Here and throughout, we will employ the useful notation
$$\minp(m,n)=\begin{cases}\min(m,n)&\text{if }m\ne n\\
\text{a number }>m&\text{if }m=n.\end{cases}$$
Note that $\minp(m,m)$ is not a well-defined number, and that $\nu(m+n)=\min'(\nu(m),\nu(n))$.
\begin{lem}\label{D2lem} $($\cite{D2}$)$ Let $\bN$ denote the set of nonegative integers. A function $f:\Z\to \bN\cup\{\infty\}$
is of the form $f(n)=\nu(n-E)$ for some $2$-adic integer $E$ iff it satisfies
$$f(n+2^d)=\minp(f(n),d)$$
for all $d\in\bN$ and all $n\in\Z$. In this case, $E=\sum_{i\ge0}2^{e_i}$, where $e_0=f(0)$ and $e_{k+1}=f(2^{e_0}+\cdots+2^{e_k})$.
\end{lem}

We begin the proofs of the theorems of Section \ref{intro} by discussing the proof of Theorem \ref{36}.
The way that the cases of Theorem \ref{36} are discovered is by having {\tt Maple} compute values of $\nu(T_n(k))$
for ranges of values of $k$. For example, seeing that $\dstyle{\nu\bigl(\sum_{i\text{ odd}}\tbinom{19}ii^k}\bigr)$ takes on the values
17, 25, 17, 18, 17, 19, 17, 18, 17, 20, 17, 18 as $k$ goes 10, 18, 26,$\ldots,98$
makes one pretty sure that for all integers $x$ we have $\nu(T_{19}(8x+2))=\nu(x-x_0)+17$ for some 2-adic integer $x_0$, and you could even
guess that
the last 9 digits in the binary expansion of $x_0$ are 100000010. But to prove it, more is required. This is a case not covered by any of our three general theorems,
but the proofs of all four of our theorems have similar structure.

Let $f(x)=\nu(T_{19}(8x+2))-17$. Then
\begin{eqnarray}&&f(x+2^d)\nonumber\\&=&\nonumber-17+\nu\bigl(\sum\tbinom{19}{2i+1}(2i+1)^{8x+2}\\
&&+\sum\tbinom{19}{2i+1}(2i+1)^{8x+2}((2i+1)^{2^{3+d}}-1)\bigr)\label{twoparts}\\
&=&\minp\bigl(f(x),\nu\bigl(\sum\tbinom{19}{2i+1}(2i+1)^{8x+2}((2i+1)^{2^{3+d}}-1)\bigr)-17\bigr).\nonumber\end{eqnarray}
Thus  the claim that $\nu(T_{19}(8x+2))=\nu(x-x_0)+17$ for some 2-adic integer $x_0$ will follow from Lemma \ref{D2lem} once we show that
\begin{equation}\label{19eq}\nu\bigl(\sum\tbinom{19}{2i+1}(2i+1)^{8x+2}((2i+1)^{2^{3+d}}-1)\bigr)=d+17\end{equation}
for all $x$ and $d\ge0$. We expand the two powers of $(2i+1)$, obtaining terms, for $k\ge0$ and $j>0$, with 2-exponent
\begin{equation}\label{term}\nu\tbinom{8x+2}k+\nu\tbinom{2^{3+d}}j+k+j+\nu\bigl(\sum\tbinom{19}{2i+1}i^{k+j}\bigr).\end{equation}
Let $\psi(s)=s+\nu(\sum\binom{19}{2i+1}i^s)$.
Since $\nu\binom{2^{3+d}}j=3+d-\nu(j)$, it will suffice to show that the minimum value of $\psi(k+j)-\nu(j)+\nu\binom{8x+2}k$ is 14,
and that this value occurs for an odd number of pairs $(k,j)$. {\tt Maple} computes that the minimum value of $\psi(s)$ is 16, which occurs
when $s=3$, 5, 7, or 9, and that $\psi(s)=17$ for $s=1$, 4, 6, 8, and 10. For $s\ge11$, $\psi(s)\ge7+s$ by \ref{DSthm}.
This information makes it easy to check that the minimum value of $\psi(k+j)-\nu(j)+\nu\binom{8x+2}k$ is indeed 14, and this value occurs
exactly when $(k,j)=(0,8)$, $(2,8)$, or $(1,8)$. This completes the proof that for all integers $x$ we have $\nu(T_{19}(8x+2))=\nu(x-x_0)+17$ for some 2-adic integer $x_0$. Each of the cases of Theorem \ref{36} can be established in this manner, although many of the
cases are covered by our general theorems \ref{thm0}, \ref{thm2}, and \ref{thm3}.

 The cases in which $T_n(k)$ is constant on a congruence class
are proved similarly, although Lemma \ref{D2lem} need not be used.
For example, to show $\nu(T_{13}(8x+7))=11$ for all $x$, we first define
$$\theta(k)=2^k\sum\tbinom{13}{2i+1}i^k.$$
{\tt Maple} and \ref{DSthm} show
$$\nu(\theta(k))\begin{cases}=10&k=5,6\\
=11&k=1,2,7\\
>11&\text{other }k.\end{cases}$$
Since $T_{13}(8x+7)=\sum\binom{8x+7}k\theta(k)$, and $\binom{8x+7}k$ is odd for $k\in\{5,6,1,2,7\}$, we obtain, mod $2^{12}$,
$$T_{13}(8x+7)\equiv\tbinom{8x+7}5\theta(5)+\tbinom{8x+7}6\theta(6)+2^{11}.$$
{\tt Maple} shows $\theta(5)\equiv\theta(6)\equiv 3\cdot2^{10}$ mod $2^{12}$. Since
$$\tbinom{8x+7}5+\tbinom{8x+7}6=\tbinom{8x+7}5(1+\tfrac{8x+2}6)\equiv0\pmod 4,$$
we obtain $\tbinom{8x+7}5\theta(5)+\tbinom{8x+7}6\theta(6)\equiv0$ mod $2^{12}$, from which our desired conclusion follows.
This concludes our comments regarding the proof of Theorem \ref{36}.

Now we work toward proofs of the more general results, Theorems \ref{thm0}, \ref{thm2}, and \ref{thm3}. First we recall some background information.
 We will often use that
\begin{equation}\label{Sa} (-1)^jj!S(k,j)=\sum\tbinom j{2i}(2i)^k-T_j(k).\end{equation}
Sometimes we have
$k<j$, in which case $S(k,j)=0$, and so $T_j(k)=\sum\binom j{2i}(2i)^k$ when $k<j$. Other times
we use (\ref{Sa}) to say that $T_j(k)\equiv\pm j!S(k,j)$ mod $2^k$.

Many times we will use without comment the
fact, related to (\ref{ale}), that
$$\nu\bigl(\tbinom mn\bigr)=\a(n)+\a(m-n)-\a(m).$$
Closely related is the fact that $\binom mn$ is odd iff each digit in the binary expansion of $m$ is at least as large as the corresponding digit of $n$.
We will sometimes say that $\binom mn$ is even due to the $2^t$-position, meaning that in this position $m$ has a 0 and $n$ has a 1.
Other basic formulas that we use without comment are
$$\a(n-1)=\a(n)-1+\nu(n)$$
and, if $0<\Delta<2^t$, then  $$\a(2^{t+1}A+2^t+\Delta)=\a(A)+t-\a(\Delta-1).$$

We also use   the well-known formula
\begin{equation}\label{Smod2}S(k+i,k)\equiv \tbinom{k+2i-1}{k-1}\mod 2.\end{equation}
A generalization to mod 4 values was given in \cite{CM} and will be used several times. We will not bother to state all eight cases
of that theorem---just those that we need.

The proof of  Theorem \ref{thm0} utilizes the following two lemmas.

\begin{lem}\label{alem} Let $1\le i\le 2^{e-1}$ with $e\ge2$. Then
$$\nu(T_{2^e+1}(2^{e-1}+i))\begin{cases}=2^e-1+i&  i\in\{ 2^{e-2},2^{e-1}\}\\
\ge 2^e+i&\text{otherwise},
\end{cases}$$
while
$$\nu(T_{2^e+2}(2^{e-1}+i))\begin{cases}=2^e-1+i& 1\le i\le 2^{e-2}\\
=2^e+i&2^{e-2}<i<2^{e-1}\\
> 2^e+i&i= 2^{e-1},
\end{cases}$$

\end{lem}

\begin{proof} For the first part, by the remarks following (\ref{Sa}), we must show
$$ \nu(\sum\tbinom{2^e+1}{2j}j^{2^{e-1}+i})\begin{cases}=2^{e-1}-1& i\in\{2^{e-2},2^{e-1}\}\\
\ge2^{e-1}&\text{otherwise},\end{cases}$$
or equivalently
$$\tfrac 1{(2^{e-1})!}\sum\tbinom{2^e+1}{2j}j^{2^{e-1}+i}\equiv\begin{cases}1\mod2& i\in\{ 2^{e-2},2^{e-1}\}\\
0\mod2&\text{otherwise}.\end{cases}$$
By Proposition \ref{stirowr}, the LHS is congruent mod 2 to $S(2^{e-1}+i,2^{e-1})$, and by (\ref{Smod2}) this is $\binom{2^{e-1}-1+2i}{2^{e-1}-1}$, which is as required.

The second part of the lemma reduces similarly to showing
\begin{equation}\label{css}{\tfrac1{(2^{e-1}+1)!}}\sum_j\tbinom{2^e+2}{2j}j^{2^{e-1}+i}\equiv\begin{cases}1\mod 2&1\le i\le 2^{e-2}\\
2\mod 4&2^{e-2}<i<2^{e-1}\\
0\mod4&i=2^{e-1}.\end{cases}\end{equation}
By \ref{stirowr}, the LHS is congruent mod 4 to \begin{equation}\label{twoter}S(2^{e-1}+i,2^{e-1}+1)+2S(2^{e-1}+i,2^{e-1}).\end{equation}
Mod 2, this is $\binom{2^{e-1}+2i-2}{2^{e-1}}$ which is odd if $1\le i\le 2^{e-2}$.
Now assume $2^{e-2}<i<2^{e-1}$.  The second term of (\ref{twoter}) is easily seen to be 0 mod 4 using (\ref{Smod2}).
For the first term of (\ref{twoter}),
we use part of \cite[Thm 3.3]{CM}, which relates mod 4 values of $S(n,k)$ to binomial coefficients. It implies that, if $e\ge3$, the mod 4 value of the first term is $\binom{2^{e-2}+k}{2^{e-3}}$, where $0\le k<2^{e-3}$.
 The 2-exponent in this number is $1+\a(2^{e-3}+k)-\a(2^{e-2}+k)=1$, as desired. If $i=2^{e-1}$, both terms of (\ref{twoter}) are 2 mod 4, by a
 similar analysis.
\end{proof}

\begin{lem}\label{Tlem} If $p\in\Z$, $\delta=1$ or $2$, and $\nu(n)=e+\Delta$ with $\Delta\ge-1$, then
$$\nu\bigl(\sum_i\tbinom{2^e+\delta}{2i+1}(2i+1)^{p}((2i+1)^n-1)\bigr)=2^e+\Delta+\delta-1.$$
\end{lem}
\begin{proof} The sum equals $\ssum_{j>0}T_j$, where
$$T_j:=2^j\tbinom{n}j\sum_i\tbinom{2^e+\delta}{2i+1}(2i+1)^{p}i^j.$$
For evaluation of the 2-exponent of the $i$-sum here, we use Corollary \ref{lcor1}. We obtain that if $j\le2^{e+\Delta}$,
then $$\nu(T_j)\ge j+e+\Delta-\nu(j)+\begin{cases}2^e+\delta-2-j&1\le j\le2^{e-1}\\
2^{e-1}-1&j>2^{e-1},\end{cases}$$
with equality if $j=2^{e-1}$.
This is $\ge 2^e+\Delta+\delta-1$ with equality iff $j=2^{e-1}$. If $j>2^{e+\Delta}$, then $\nu(T_j)>2^e+\Delta$ since
$2^{e+\Delta}+2^{e-1}>2^e+\Delta$ for $\Delta\ge-1$.
\end{proof}

Now we easily prove Theorem \ref{thm0}.
\begin{proof}[Proof of Theorem \ref{thm0}] Let $\delta\in\{1,2\}$, $1\le i\le 2^{e-1}$, and let
$$g(x)=\nu(T_{2^e+\delta}(2^{e-1}x+2^{e-1}+i))-2^e+2-\delta.$$
Note that the expression that we wish to evaluate for Theorem \ref{thm0} is $g(x-1)+2^e-2+\delta$.

For $d\ge0$, writing $T_n(-)$ as a sum of two parts as we did in (\ref{twoparts}),
$$g(x+2^d)=\min\nolimits'\bigl(g(x),-2^e+2-\delta+\nu(\sum\tbinom{2^e+\delta}{2j+1}(2j+1)^p((2j+1)^{2^{d+e-1}}-1))\bigr),$$
where $p=2^{e-1}x+2^{e-1}+i$. By Lemma \ref{Tlem}, the RHS equals $\min'(g(x),d)$, and so $g(x)=\nu(x-E_\delta)$ for some $E_\delta$ by Lemma \ref{D2lem}. By Lemma \ref{alem}
$$\nu(E_\delta)=g(0)\begin{cases}=i&\delta=1,\ i\in\{2^{e-2},2^{e-1}\}\\
>i&\delta=1,\ i\not\in\{2^{e-2},2^{e-1}\}\\
=i-1&\delta=2,\ 1\le i\le2^{e-2}\\
=i&\delta=2,\ 2^{e-2}<i<2^{e-1}\\
>i&\delta=2,\ i=2^{e-1}.\end{cases}$$
Our desired $g(x-1)+2^e-2+\delta$ equals $\nu(x-1-E_\delta)+2^e-2+\delta$, and $x_{i,2^e+\delta}:=1+E_\delta$ is as claimed.
\end{proof}

The proof of Theorem \ref{thm2} is similar in nature, but longer.
\begin{proof}[Proof of Theorem \ref{thm2}] Using Lemma \ref{D2lem} and arguing as in (\ref{19eq}), it suffices to prove that for $d\ge0$ and any integer $x$
\begin{equation}\label{abeq}\nu\bigl(\sum_i\tbinom n{2i+1}(2i+1)^{2^{e+a}x}((2i+1)^{2^{e+a+d}}-1)\bigr)=n-b+d.\end{equation}
Indeed, if
$$g(x)=\nu\bigl(\sum_i\tbinom n{2i+1}(2i+1)^{2^{e+a}x}\bigr)-n+b,$$
then (\ref{abeq}) implies $g(x+2^d)=\min'(g(x),d)$  and Theorem \ref{thm2} then follows
from Lemma \ref{D2lem}.

We write the sum in (\ref{abeq}) as $\sum T_j$ with
\begin{equation}\label{nut}\nu(T_j)=j+e+a+d-\nu(j)+\nu\bigl(\sum_i\tbinom n{2i+1}(2i+1)^{2^{e+a}x}i^j\bigr).\end{equation}
We will show that in all cases $\nu(T_j)$ is minimized for a unique value of $j$.

The {\bf second case} of the theorem will follow from proving that if $2^{e}<n\le3\cdot2^{e-1}$ and $\nu(p)\ge e-1$, then
$$j+e-1-\nu(j)+\nu\bigl(\sum\tbinom n{2i+1}(2i+1)^pi^j\bigr)\ge n-2$$
with equality iff $j=2^{e-1}$. Expanding $(2i+1)^p$ as $\sum_{k\ge0} 2^k\binom pk i^k$ leads us to needing that for $j>0$
\begin{equation}\label{cs1a}j+e-\nu(j)+\nu\bigl(\sum\tbinom n{2i+1}i^j\bigr)\ge n-1\end{equation}
with equality iff $j=2^{e-1}$, and
\begin{equation}\label{cs1}j+k+2e-\nu(j)-\nu(k)+\nu\bigl(\sum\tbinom n{2i+1}i^{j+k}\bigr)> n\end{equation}
for $j,k>0$.

The equality in (\ref{cs1a}) when $j=2^{e-1}$ follows easily from Proposition \ref{qprop2}.
Also by \ref{qprop2}, the difference in (\ref{cs1a}) becomes
\begin{eqnarray}&&j+e+1-\nu(j)+\nu\bigl(\sum\tbinom n{2i+1}i^j\bigr)-n\label{jee}\\
&\ge&e-\a(j)-\nu(j)=e-1-\a(j-1)\nonumber.\end{eqnarray}
This is $>0$ if $j\ne2^{e-1}$ and $j<3\cdot2^{e-2}$, while if $j=3\cdot2^{e-2}$, then $\binom{n-1-j}j=0$
and so (\ref{jee}) is $>0$ by \ref{qprop2}.

Now suppose $j>3\cdot2^{e-2}$. Then $j-\nu(j)>3\cdot2^{e-2}$, and since $n\le 3\cdot2^{e-1}$,
\begin{equation}\label{47}n-\nu([\tfrac n2]!)\le3\cdot2^{e-2}+e-1.\end{equation}
Thus, using Proposition \ref{DSthm}, we obtain
$$j+e-1+\nu(j)+\nu\bigl(\sum\tbinom n{2i+1}i^j\bigr)>3\cdot2^{e-2}+e+1+\nu([\tfrac n2]!)\ge n,$$
establishing strict inequality in (\ref{cs1a}).

Now we verify (\ref{cs1}). By \ref{qprop2}, (\ref{cs1}) is satisfied if
\begin{equation}\label{dn1}\nu(j)+\nu(k)+\a(j+k)\le 2e-2\end{equation}
or if
\begin{equation}\label{dn2}\nu(j)+\nu(k)+\a(j+k)=2e-1\text{ and }\tbinom{n-1-j-k}{j+k}\equiv0\pmod2.\end{equation}
By \ref{DSthm} and (\ref{47}), (\ref{cs1}) is also satisfied if
\begin{equation}\label{dn3}j+k-\nu(j)-\nu(k)\ge3\cdot2^{e-2}-e.\end{equation}

If $j+k>3\cdot2^{e-2}$, then (\ref{dn3}) is satisfied. If $\{j,k\}=\{2^{e-1},2^{e-2}\}$, then (\ref{dn2}) is satisfied.
Assume WLOG that $\nu(j)\ge\nu(k)$. Then (\ref{dn1}) is implied by $\nu(j)+1+\a(j+k-1)\le 2e-2$ and this is satisfied
whenever $j+k\le3\cdot2^{e-2}$ and $(j,k)\ne(2^{e-1},2^{e-2})$.

The {\bf third case} of the theorem will follow from proving that, referring to (\ref{nut}), if $2^{e+1}-2^{t+1}<n\le 2^{e+1}-2^t$ with $1\le t<e-1$
and $\nu(p)\ge e$, then
$$j+e+d-\nu(j)+\nu\bigl(\sum\tbinom n{2i+1}(2i+1)^pi^j\bigr)\ge n-1+d$$
with equality iff $j=2^{e}-2^t$. Expanding $(2i+1)^p$, this reduces to showing if $j>0$ then
\begin{equation}\label{one}j+e+1-\nu(j)+\nu\bigl(\sum\tbinom n{2i+1}i^j\bigr)\ge n\end{equation}
with equality iff $j=2^{e}-2^t$, and if $j,k>0$, then
\begin{equation}\label{two}j+k+2e-\nu(j)-\nu(k)+\nu\bigl(\sum\tbinom n{2i+1}i^{j+k}\bigr)\ge n.\end{equation}

If $j>2^{e}$, since $n\le 2^{e+1}-2$,
$$j+e+1-\nu(j)\ge2^{e}+e+1>n-\nu([\tfrac n2]!),$$
and so strict inequality holds in (\ref{one}) by \ref{DSthm}.

By Theorem \ref{qprop2}, (\ref{one}) is satisfied if
\begin{equation}\label{onep}e\ge\nu(j)+\a(j)=\a(j-1)+1,\end{equation}
and equality holds in (\ref{one}) iff equality holds in (\ref{onep}) and $\binom{n-1-j}j$ is odd.
If $j<2^{e}$, then $\a(j-1)\le e-1$ with equality iff $j=2^{e}-2^r$ for some $r$. Thus (\ref{onep})
holds with equality iff $j=2^{e}-2^t$ by Lemma \ref{bcl}.

If $j=2^{e}$, by Proposition \ref{refine} the LHS of (\ref{one}) is $\ge n+1$.
Thus (\ref{one}), including consideration of equality, has been established for all $j$.

By \ref{DSthm}, (\ref{two}) is satisfied if
$$j+k+2e-\nu(j)-\nu(k)\ge n-\nu([\tfrac n2]!),$$
and hence, since $n\le 2^{e}-2$, it is satisfied if
$$j+k+2e-\nu(j)-\nu(k)\ge 2^{e}+e-1.$$
This is satisfied if $j+k>2^{e}$.

By \ref{qprop2}, (\ref{two}) is also satisfied if
$$\nu(j)+\nu(k)+\a(j+k)\le 2e-1.$$
This is satisfied if $j=k=2^{e-1}$ and if $\nu(j)+\a(j+k-1)\le 2e-2$, which is true for all other $(j,k)$ with $j+k\le 2^{e}$.

The {\bf first case}, $n=2^{e}$, will follow similarly from
\begin{equation}\label{j}j+e-1-\nu(j)+\nu\bigl(\sum\tbinom{2^{e}}{2i+1}i^j\bigr)\ge 2^{e}\end{equation}
for $j>0$ with equality iff $j=2^{e-2}$, while if $j,k>0$, then
\begin{equation}\label{jk}j+k+2e-2-\nu(j)-\nu(k)+\nu\bigl(\sum\tbinom{2^{e}}{2i+1}i^{j+k}\bigr)\ge2^{e}+2.\end{equation}

Equality in (\ref{j}) with $j=2^{e-2}$ follows from Proposition \ref{refine} since $\binom{2^{e}-1-2^{e-2}}{2^{e-2}}\equiv 2$ mod 4.
If $j>2^{e-1}$, then strict inequality in (\ref{j}) is implied by \ref{DSthm}. If $j=2^{e-1}$, it is implied by Proposition \ref{stirowr}.
It is implied by Theorem \ref{qprop1} if $\nu(j)\le e-3$, which is true for $j<2^{e-1}$ provided $j\ne2^{e-2}$.

Similarly, (\ref{jk}) is implied by \ref{DSthm} if $j+k\ge 2^{e-1}$ unless $j=k=2^{e-2}$, in which case it is implied by \ref{stirowr}.
If $j+k<2^{e-1}$, then $\nu(j)+\nu(k)\le 2e-5$, and so the claim follows from Proposition \ref{qprop1}.
\end{proof}

The following lemma was used in the above proof.
\begin{lem}\label{bcl} If $2^{e+1}-2^{t+1}\le m<2^{e+1}-2^t$ with $0\le t<e$ and $j=2^{e}-2^r$ with $0\le r<e$,
then $\binom{m-j}j$ is odd iff $r=t$.\end{lem}
\begin{proof} If $r<t$, then $0\le m-j<j$, so $\binom{m-j}j=0$. If $r=t$, then $\binom{m-j}j=\binom{2^{e}-2^t+d}{2^{e}-2^t}$ with $0\le d<2^t$
and hence is odd. If $r>t$, then the binary expansion of $m-j$ has a 0 in the $2^r$ position, while $j$ has a 1 there.\end{proof}

The following lemma will be useful in the proof of Theorem \ref{thm3}.

\begin{lem}\label{plem} In the notation of Theorem \ref{thm3}, if $p\in S_n$, then $\a(p-p_0)\le1$
and $\binom{n-p_0-1}{p_0}$ and $\binom{n-2^{e-1}-p_0-1}{2^{e-1}+p_0}$ are odd.\end{lem}
\begin{proof} Let $p=A2^t+p_0$ and $n=2^e+2^t+\Delta$ with both $p_0$ and $\Delta$ nonnegative and less than $2^t$. If $\Delta-p_0-1<0$, then $\binom{2^e-A2^t}{A2^t}$ is odd,
which easily implies $\a(A)\le1$, while if $\Delta-p_0-1\ge0$, then $\binom{2^e+2^t-A2^t}{A2^t}$ is odd. This implies that $A$ is 0 or an even 2-power.

Now, if $p\ne p_0$, we can write $p=p_0+2^{t+r}$ with $r\ge0$. If $r=0$, then $\binom{2^e+\Delta-1-p_0}{2^t+p_0}$ odd implies $\Delta-1-p_0<0$
and $\binom{2^t+\Delta-1-p_0}{p_0}$ odd, which implies $\binom{n-p_0-1}{p_0}$ odd. If $r>0$, then the odd binomial coefficient can be considered
mod 2 as $\binom{2^e-2^{t+r}}{2^{t+r}}\binom{2^t+\Delta-1-p_0}{p_0}$, which implies $\binom{n-p_0-1}{p_0}$ is odd.

Now we may assume $\binom{n-p_0-1}{p_0}$ is odd. Thus $\binom{2^t+\Delta-1-p_0}{p_0}$ is odd and hence so is $\binom{2^{e-1}+2^t+\Delta-1-p_0}{2^{e-1}+p_0}$.
\end{proof}

The proof of Theorem \ref{thm3} is similar to the others, but longer yet.
\begin{proof}[Proof of Theorem \ref{thm3}]
Similarly to the proofs of the other three theorems, it suffices to prove for $d\ge0$ and any integer $x$
\begin{equation}\label{t3e}\nu\bigl(\sum\tbinom n{2i+1}(2i+1)^{2^{e-1}x+q}((2i+1)^{2^{e-1+d}}-1)\bigr)=d+n-2-\a(p_0).\end{equation}
Here, and for the remainder of this section, $n$, $e$, $t$, $q$, $p$, and $p_0$ are as in Theorem \ref{thm3}. To prove (\ref{t3e}), it suffices to show for $k\ge0$ and $j>0$
\begin{equation}\label{t3e2}\nu\tbinom{2^{e-1}x+q}k+e-1-\nu(j)+j+k+\nu\bigl(\sum\tbinom n{2i+1}i^{j+k}\bigr)\ge n-2-\a(p_0)\end{equation}
with equality iff $j=2^{e-1}$ and $k=p_0$.

We first prove the equality. Note that if $q$ is associated to $p\in S_n$, then $\binom{2^{e-1}x +q}{p_0}$ is odd. We must show that
$$\nu\bigl(\sum\tbinom n{2i+1}i^{2^{e-1}+p_0}\bigr)=n-2-p_0-2^{e-1}-\a(p_0).$$
This follows from Proposition \ref{qprop2} since $\binom{n-2^{e-1}-1-p_0}{2^{e-1}+p_0}$ is odd by Lemma \ref{plem}.


Strict inequality in (\ref{t3e2}) when $j=2^{e-1}$ and $k\ne p_0 $ follows from Lemma \ref{techlem} using Propositions
\ref{DSthm} and \ref{refine}. Here we also use that if $p\in S_n$ and $k$ satisfies (\ref{kres}) then $k<2^{e-1}$ and hence the $x$ in
$\binom{2^{e-1}x+q}k$ does not play an essential role.
\begin{lem}\label{techlem} If $n$, $e$, $t$, $q$, $p$, and $p_0$ are as in Theorem \ref{thm3} and
\begin{equation}\label{kres}0\le k\le n-\nu([\tfrac n2]!)-2-2^{e-1}-\a(p_0),\end{equation}
then
\begin{equation}\label{alpheq}\a(q-k)+\a(p_0)-\a(q)\ge-1.\end{equation}
If the LHS of (\ref{alpheq}) equals $-1$, then either $n$ is even and $\binom{n-1-2^{e-1}-k}{2^{e-1}+k}\equiv0$ mod $4$ or $\binom{n-1-2^{e-1}-k}{2^{e-1}+k}=0=\binom{n-1-2^{e-1}-k}{2^{e-1}+k-2}$.
If the LHS of (\ref{alpheq}) equals $0$, then either $\binom{n-1-2^{e-1}-k}{2^{e-1}+k}\equiv0$ mod $2$ or $k=p_0$.
\end{lem}
\begin{proof} We begin by proving (\ref{alpheq}).   Using Lemma \ref{plem} and that $q$ is associated to $p$, we have
\begin{equation}\label{newqdef}q=p_0+\delta2^{t+r}+\eps2^w\end{equation}
with $\delta$ and $\eps$ equal to 0 or 1, $r\ge0$, and $w=\nu(n)-1$ or $w>t$. The only way that (\ref{alpheq}) could fail is
if $k=q=p_0+2^{t+r}+2^w$. But (\ref{kres}) implies $k\le 2^t+t-2$, which is inconsistent with $w>t$ and with $r>0$. Thus $n$ is even and
$w=\nu(n)-1$ and $r=0$. Let $n=2^e+2^t+c2^{t-1}+2b$ with $c\in\{0,1\}$ and $b\le 2^{t-2}-1$. If $c=0$, then (\ref{kres})
reduces to $p_0+\a(p_0)+2^{t-2}+2^w \le t-3$, which is false. If, on the other hand, $c=1$, the assumption that $\binom{n-p-1}p$ is odd and $p\ge2^t$ implies
that $p_0\ge2^{t-1}+2b$, so $k>2^t+2^{t-1}$ contradicts $k\le 2^t+t-2$.

There are three conceivable ways in which equality could hold in (\ref{alpheq}). One is $\eps=0$, $\delta=1$, and $k=q$; i.e., $k=q=p\ge2^t$.
But $k\ge2^t$ implies $2^{e-1}+k>n-1-2^{e-1}-k$ and hence $\binom{n-1-2^{e-1}-k}{2^{e-1}+k}=0$. We also have $\binom{n-1-2^{e-1}-k}{2^{e-1}+k-2}=0$;
the only way this could fail is if $n=2^e+2^{t+1}-1$ and $k=p=2^t$, but then $\binom{n-p-1}p$ is not odd.
Another is $\eps=\delta=1$ and $\a(q-k)=1$. In this case, the only way to have $k<2^t$ is if $w=\nu(n)-1$ and $k=q-2^{t+r}$, where $r$ is as in (\ref{newqdef}). In this case,
$\binom{n-1-2^{e-1}-k}{2^{e-1}+k}\equiv0$ mod 4, using the result that $\binom{a+b}b$ is divisible by $2^t$ if there are at least $t$
carries in the binary addition of $a$ and $b$. In this case, either the binomial coefficient equals 0, or else, if $v=\nu(n)$, there will be carries in the $2^{v-1}$ and $2^v$ positions in the
relevant binary addition.  The third possibility, $\eps=1$, $\delta=0$, and $k=q$, implies that
 $n$ is even and $\nu(q)=\nu(n)-1$ and leads to $\binom{n-1-2^{e-1}-k}{2^{e-1}+k}\equiv0$ mod 4, exactly as above.

Finally we show that if $\binom{n-1-2^{e-1}-k}{2^{e-1}+k}$ is odd and the LHS of (\ref{alpheq}) equals 0, then $k=p_0$.
It is not difficult to see that if $0\le k<n-2^{e-1}$ and $\binom{n-1-2^{e-1}-k}{2^{e-1}+k}$ is odd, then $k<2^t$ and $\binom{n-1-k}k$ is odd.

First suppose $\a(q-k)=2$ and $\a(q)-\a(p_0)=2$. If $q=p_0+2^{t+r}+2^{t+s}$ with $0\le r<s$, then to keep $k<2^t$, we must have $k=p_0$.
If $q=p_0+2^{t+r}+2^{\nu(n)-1}$ with $r\ge0$, then we must have $k=p_0+2^{\nu(n)-1}-2^u$ for some $u$. If $u\ne \nu(n)-1$, then $\binom{n-1-2^{e-1}-k}{2^{e-1}+k}$ is even due to the $2^{\min(u,\nu(n)-1)}$-position.

Now suppose $\a(q-k)=1$ and $\a(q)-\a(p_0)=1$. If $q=p_0+2^{t+r}$ with $r\ge0$, then $k=p_0$ is the only way to have $k<2^t$.
If $q=p_0+2^{\nu(n)-1}$, then the argument at the end of the preceding paragraph applies. This completes the proof of Lemma \ref{techlem},
and hence the proof that when $j=2^{e-1}$, (\ref{t3e2}) holds with equality exactly as claimed there.
\end{proof}

We continue the proof of Theorem \ref{thm3} by establishing strict inequality in (\ref{t3e2}) when $0<j<2^{e-1}$ and $0\le k<2^{e-1}$.
The following elementary lemma will be useful.
\begin{lem}\label{jklem} Suppose $0<j<2^{e-1}$ and $0\le k<2^{e-1}$. Let $\phi(j,k)=\a(j+k)+\nu(j)-\a(k)$.
Then
\begin{enumerate}
\item $\phi(j,k)\le e-1$;
\item $\phi(j,k)=e-1$ iff $j=2^{e-1}-2^h$ and $0\le k<2^h$ for some $0\le h<e-1$;
\item $\phi(j,k)=e-2$ iff either $j=2^{e-1}-2^h$ and $2^{e-2}\le k<2^{e-2}+2^h$ for some $0\le h<e-1$, or $j=2^{e-1}-2^\ell-2^h$, $0\le h\le \ell<e-1$, and $0\le k<2^h$ or $2^\ell\le k<2^\ell+2^h$.
    \end{enumerate} \end{lem}
\begin{proof} Let $h=\nu(j)$, and let $k_0=k-(k\mod 2^h)$. Then $\phi(j,k)=h+\a(j+k_0)-\a(k_0)$.
The only way to get $\a(j+k_0)-\a(k_0)=e-h-1$  is if $j+k_0=2^e-2^h$ and $\a(k_0)=1$, or $j+k_0=2^{e-1}-2^h$ and $k_0=0$. But the first is impossible since $k_0<2^{e-1}$. Similarly the only ways to get $\a(j+k_0)-\a(k_0)=e-h-2$ is if $(\a(j+k_0),\a(k_0))=(e-h-1,1)$ or $(e-h,2)$, and these
can only be accomplished in the ways listed in part (3).
\end{proof}

Let
$${\tbinom{m-s}s}'=\tbinom{m-s}s+\begin{cases}2\tbinom{m-s}{s-2}&\text{if $m$ and $s$ are even}\\
0&\text{otherwise}\end{cases}$$
and \begin{equation}\label{nutdef}\nut_2\tbinom{m-s}s=\min(2,\nu({\tbinom{m-s}s}')).\end{equation}
Let $\phi(j,k)$ be as in Lemma \ref{jklem}. The desired strict inequality in (\ref{t3e2}) when $0<j<2^{e-1}$ and $0\le k<2^{e-1}$
follows from the following result using Proposition \ref{refine}.
\begin{thm}\label{e-1thm} If $n,e,t,q$, $p$, and $p_0$ are as in Theorem \ref{thm3}, $0<j<2^{e-1}$, and $0\le k<2^{e-1}$,  then
\begin{equation}\label{impineq}\a(q-k)+e-1+\nut_2\tbinom{n-1-j-k}{j+k}\ge\a(q)-\a(p_0)+\phi(j,k).\end{equation}
\end{thm}
\begin{proof} By Lemma \ref{plem}, $p=p_0$ or $p_0+2^{t+s}$ with $s\ge0$. Hence $\a(q)-\a(p_0)\le2$. Also $\nu(p)\ge\nu(n)$, a consequence of the oddness of $\tbinom{n-1-p}p$, will be used often without comment. The theorem will follow from showing:
\begin{itemize}
\item if $\phi(j,k)=e-1$, then
$$\nut_2\tbinom{n-1-j-k}{j+k}\ge\begin{cases}2&\text{if $\a(q)-\a(p_0)=2$ and $k=q$}\\
1&\text{if $\a(q)-\a(p_0)=2$ and $\a(q-k)=1$}\\
1&\text{if $\a(q)-\a(p_0)=1$ and $k=q$,}\end{cases}$$
\item and if $\phi(j,k)=e-2$, then
$$\nu\tbinom{n-1-j-k}{j+k}\ge1\text{ if $\a(q)-\a(p_0)=2$ and $k=q$.}$$
\end{itemize}
We call these cases 1 through 4. Let $n=2^e+2^t+\Delta$ with $0\le\Delta<2^t$. Our hypothesis is that
$\binom{2^e+2^t-\eps2^{t+s}+\Delta-p_0-1}{\eps2^{t+s}+p_0}$ is odd.

{\bf Case 3:} We have $q=p_0+2^{r}$ with $r\ge t$ or  $r=\nu(n)-1$, in which latter case $\Delta$ and $p_0$ are divisible by $2^{r+1}$.
We must show that $\binom{2^{e-1}+2^t+2^h+\Delta-1-p_0-2^r}{2^{e-1}-2^h+p_0+2^r}$ is even. Here $2^h>p_0+2^r$.
If $r\ge t$, then the binomial coefficient is even due to the $2^r$- or $2^h$-position, while if $r=\nu(n)-1$, it is even due to the $2^{\nu(n)-1}$-position.

{\bf Case 2:} Here $q=p_0+2^s+2^r$ with $s\ge t$ and $r=\nu(n)-1$ or $r>s$. Also $k=q-2^v$ and $2^h>k$. The binomial coefficient which we must show is even is
$$C:=\binom{2^{e-1}+2^h+2^t+\Delta-1-p_0-2^s-2^r+2^v}{2^{e-1}-2^h+p_0+2^s+2^r-2^v}.$$
If $v=r$ or $s$, it reduces to Case 3, just considered. If $r=\nu(n)-1$, then $C$ is even due to the $2^{\min(v,\nu(n)-1)}$-position.
Otherwise $C$ is even due to the $2^h$-position, since $2^t+\Delta-1-p_0-2^s-2^r+2^v$ is negative.

{\bf Case 1:} Now $q$ is as in Case 2, but $k=q$. We must show that there are at least two carries in the binary addition of
$2^{e-1}-2^h+p_0+2^s+2^r$ and $2^{h+1}+2^t+\Delta-1-2p_0-2^{s+1}-2^{r+1}$. If $r=\nu(n)-1$, carries occur in positions $2^r$ and $2^{r+1}$.
If $r>s$, carries occur in $2^r$ and $2^h$. The second term in the definition of ${\binom{m-s}s}'$ is easily seen to be inconsequential here.

{\bf Case 4:} Again $q$ is as in Case 2, $k=q$, and $(j,k)$ is one of the two types in Theorem \ref{e-1thm}. For the first type of $(j,k)$, if $r>s$, then $j+k>n-1-j-k$, so $\binom{n-1-j-k}{j+k}=0$, while if $r=\nu(n)-1$, the binomial coefficient is even due to the $2^{\nu(n)-1}$-position.
If $(j,k)$ is of the second type and $k=p_0+2^s+2^{\nu(n)-1}$, then $\binom{n-1-j-k}{j+k}$ is even due to the $2^{\nu(n)-1}$-position,
since $p_0$, $n$, and $j$ are all divisible by $2^{\nu(n)}$.

If $j=2^{e-1}-2^\ell-2^h$ and $2^\ell\le k<2^\ell+2^h$ with $k=p_0+2^s+2^r$ with $r>s$, we claim that $\binom{n-1-j-k}{j+k}$ is even due to the $2^{e-2}$-position. Indeed, $2^{e-1}-2^h\le j+k<2^{e-1}$, so $j+k$ has a 1 in the $2^{e-2}$-position, while $$2^{e-1}\le n-j-k-1<2^{e-1}+2^{t+1}-2^h<2^{e-1}+2^{e-2}$$
since $t<e-3$. If $j=2^{e-1}-2^{e-2}-2^h$ and $2^{t+1}< k<2^h$, one easily verifies that $\binom{n-1-j-k}{j+k}$ is even due to the $2^{e-3}$-position. Finally, if $j=2^{e-1}-2^\ell-2^h$ with $h<\ell<e-2$ and $2^{t+1}<k<2^h$, then $\binom{n-1-j-k}{j+k}$ is even due to the $2^{e-2}$-position,
as is easily proved.
\end{proof}

Our final step in the proof of Theorem \ref{thm3} is to prove strict inequality in (\ref{t3e2}) when $j>2^{e-1}$. Proposition \ref{DSthm}
implies the result if $k\ge2^t$ or if $j> 2^e$. Thus, by Proposition \ref{refine}, it suffices to prove
(\ref{impineq}) when $2^{e-1}<j\le2^{e}$ and $0\le k<2^t$. Recall that $q$ is as in (\ref{newqdef}).
Because $k<2^t$, it must be the case that if $\delta=1$, then $2^{t+r}$ appears in $q-k$, and similarly $2^w$ if $\eps=1$ and $w>t$.
These will contribute to $\a(q-k)$. Thus the only ways to have
$D_k:=\a(q-k)-(\a(q)-\a(p_0))\le0$ are (a) $k=p_0$ and $D_k=0$; (b) $k=p_0+2^{\nu(n)-1}$ and $D_k=-1$; and (c) $k=p_0+2^{\nu(n)-1}-2^v$ and
$D_k=0$.

Similarly to Lemma \ref{jklem}, we have for $2^{e-1}\le j\le2^e$ and $0\le k<2^{e-1}$, $\phi(j,k)\le e$ with equality iff
$j=2^e-2^h$ and $0\le k<2^h$ for some $0\le h<e$, or $j=2^e$.
We will be done once we prove the following lemma.
\end{proof}
\begin{lem} If $p\in S_n$ and $2^{e-1}<j\le 2^e$ and $0\le k<2^t$, then
\begin{enumerate}
\item if $k=p_0$ or $p_0+2^{\nu(n)-1}-2^v$ for some $v$, and $\phi(j,k)=e$, then $\binom{n-1-j-k}{j+k}$ is even.
\item if $k=p_0+2^{\nu(n)-1}$ and $\phi(j,k)=e$, then $\nut_2\binom{n-1-j-k}{j+k}=2$.
\item if $\phi(j,k)=e-1$ and $k=p_0+2^{\nu(n)-1}$, then $\binom{n-1-j-k}{j+k}$ is even.
\end{enumerate}
\end{lem}
\begin{proof} If $j=2^e$, then $\nut\binom{n-1-j-k}{j+k}=0=\nu\binom{n-1-j-k}{j+k}$, so now we may assume $j<2^e$.
 If $k=p_0$ or $p_0+2^{\nu(n)-1}$ and $\phi(j,k)=e$, then $j+k>n-1-j-k$ (and hence $\binom{n-1-j-k}{j+k}=0$) unless $h=e-2$ and $t=e-1$.
But part of the definition of $S_n$ said that if $t=e-1$, then $p_0\ge\Delta$, and hence $j+k>n-1-j-k$ in this case, too.
For part (2), we also need that $\binom{n-1-j-k}{j+k-2}$ is even, but it will also be 0, using that $2^{\nu(n)}-2\ge0$.

If $k=p_0+2^{\nu(n)-1}-2^v$, then $\binom{n-1-j-k}{j+k}$ is even due to the $2^{\min(v,\nu(n)-1)}$-position if $v\ne \nu(n)-1$,
while if $v=\nu(n)-1$, we are in the case $k=p_0$ already handled. A similar argument works for part (3), using the
$2^{\min(\nu(j),\nu(k))}$-position, provided $\nu(j)\ne\nu(k)$. However, equality of $\nu(j)$ and $\nu(k)$ will not occur, because one can easily prove by
induction on $j$ that if $2^{e-1}\le j<2^e$ and $0\le k<2^{e-1}$ and $\nu(j)=\nu(k)$, then $\a(j+k)+\nu(j)-\a(k)<e-1$.
\end{proof}

\section{Proofs of results about $\nu\bigl(\sum\binom n{2i+1}i^k\bigr)$}\label{sumpfsec}
In this section, we prove four propositions about $\nu\bigl(\sum\binom n{2i+1}i^k\bigr)$ which were stated and used in the previous
section. The polynomials $q_m(x)$ which we introduce in Definition \ref{qdef} might be of independent interest.

Our first proof utilizes an argument of Sun (\cite{ZW}).
\begin{proof}[Proof of Proposition \ref{stirowr}] We mimic the argument in the proof of \cite[Thm 1.3]{ZW}. Let $C_{m,\ell,b}={\displaystyle\sum\limits_{i}\tbinom m{2i+b}\tbinom i\ell}$.
Using an identity which relates $i^k$ to Stirling numbers, we obtain
\begin{eqnarray*} \sum_i\tbinom{2n+\eps}{2i+b}i^k&=&\sum_i\tbinom{2n+\eps}{2i+b}\sum_\ell\tbinom i{\ell}\ell!S(k,\ell)\\
&=&\sum_\ell C_{2n+\eps,\ell,b} \ell! S(k,\ell).\end{eqnarray*}
Since $C_{2n+1,n,0}=2n+1$, $C_{2n+1,n-1,0}=\frac23(2n+1)(n+1)n$, $C_{2n,n,0}=1$,  $C_{2n,n-1,0}=2n^2$, $C_{2n,n,1}=0$, $C_{2n,n-1,1}=2n$,
$C_{2n+1,n,1}=1$, and $C_{2n+1,n-1,1}=2n(n+1)$, our result follows from
$$\nu(\ell!C_{2n+\eps,\ell,b})\ge\nu((2n+\eps)!)-\ell=\nu(n!)+n-\ell,$$
where we have used \cite[Thm 1.1]{DS2} at the first step.
\end{proof}

The remaining proofs utilize a new family of polynomials $q_m(x)$.
\begin{defin}\label{qdef}  For $m\ge1$, we define polynomials $q_m(x)$ inductively by $q_1(x)=x-1$, and
\begin{equation}\label{q1d}\text{if }(x+1)x(x-1)\cdots(x-m+2)=\sum_{j=1}^m b_{j,m}x^j,\end{equation}
\begin{equation}\label{q2d}\text{then }(x+1)x(x-1)\cdots(x-m+2)=\sum_{j=1}^m 2^{m-j}b_{j,m}q_j(x).\end{equation}
\end{defin}
\noindent For example, $q_2(x)=x^2-x+2$. The relevance of these polynomials is given by the following result.

\begin{thm}\label{qthm2}  For all integers $x$,
$$\sum_i i^m\tbinom {x+1}{2i+1}=2^{x-2m}q_m(x).$$
\end{thm}
\begin{proof}  The proof is by induction on $m$. Validity when $m=1$ follows from
$$2\sum i\tbinom{x+1}{2i+1}+2^x=\sum(2i+1)\tbinom{x+1}{2i+1}=(x+1)\sum\tbinom x{2i}=(x+1)2^{x-1}.$$

We show that $2^{2m-x}\sum i^m\binom{x+1}{2i+1}$ satisfies the equation (\ref{q2d}) which defines $q_m(x)$.
We insert this expression for $q_j(x)$ into the RHS of (\ref{q2d}) and obtain
\begin{eqnarray*}&&2^{m-x}\sum_i\tbinom{x+1}{2i+1}\sum_j(2i)^jb_{j,m}\\
&=&2^{m-x}\sum\tbinom{x+1}{2i+1}(2i+1)\cdots(2i-m+2)\\
&=&2^{m-x}(x+1)\cdots(x-m+2)\sum\tbinom{x-m+1}{x-2i},\end{eqnarray*}
but $\sum\tbinom{x-m+1}{x-2i}=2^{x-m}$, since it is the sum of all $\tbinom{x-m+1}j$ with $j$ in a fixed parity.
Thus we obtain $(x+1)\cdots(x-m+2)$, as desired. At the second step above, we have used (\ref{q1d}) with $x=2i$.
\end{proof}

Proposition \ref{qprop1} is an immediate consequence of Theorems \ref{qthm2} and \ref{q1e}.
\begin{thm}\label{q1e} For all positive integers $x$, $$\nu(q_m(x))\ge m-x+\nu((x+1)!)=m+1-\a(x+1).$$
\end{thm}
\begin{proof}
 The proof is by induction on $m$. When $m=1$, it reduces to $\a(x+1)+\nu(x-1)\ge2$.

For the LHS of (\ref{q2d}), note that
$$\nu((x+1)\cdots(x-m+2))\ge\nu((x+1)!)-(x-m),$$
using (\ref{ale}). For the $j$-term ($j<m$) in the sum in (\ref{q2d}), by induction on $m$ we have 2-exponent
$$\ge m-j+\nu(b_{j,m})+j-x+\nu((x+1)!)\ge m-x+\nu((x+1)!).$$
Thus the inequality for $\nu(q_m(x))$ follows by induction.
 \end{proof}

The proof of Proposition \ref{qprop2} requires the following two lemmas, and the result follows easily from the second and Theorem \ref{qthm2}.
\begin{lem}\label{blem} If $b_{j,m}$ is as in Definition \ref{qdef}, then
$$\nu(b_{j,m})\ge\nu(m!)-\nu(j!)-(m-j)$$
with equality iff $\binom j{m-j}$ is odd.\end{lem}
\begin{proof} We have
\begin{eqnarray*}&&\sum_{j\ge0}x^j\sum_{m\ge j}b_{j,m}\tfrac{z^m}{m!}\\
&=&\sum_{m\ge0}{\tfrac{z^m}{m!}}\sum_{j=0}^m b_{j,m}x^j\\
&=&\sum_{m\ge0}{\tfrac 1{m!}}(x+1)_m z^m=\sum_{m\ge0}\tbinom{x+1}mz^m\\
&=&(1+z)^{x+1}=e^{(x+1)\log(1+z)}\\
&=&\sum_{k\ge0}\tfrac 1{k!}(\log(1+z))^k(1+x)^k\\
&=&\sum_{k\ge0}{\tfrac 1{k!}}(\log(1+z))^k\sum_{i=0}^k\tbinom ki x^i.\end{eqnarray*}
Here we have introduced the notation $(x+1)_m=(x+1)x\cdots(x-m+2)$.
Equate coefficients of $x^jz^m$, and get
$${\tfrac 1{m!}}b_{j,m}={\tfrac 1{j!}}\sum_k\tfrac 1{(k-j)!}([z^m](\log(1+z))^k).$$
Here $[z^m]p(z)$ denotes the coefficient of $z^m$ in $p(z)$.
Let $\ell(z)=\log(1+z)/z$. The claim of the lemma reduces to
$$\nu\bigl(\sum_k\tfrac 1{(k-j)!}([z^{m-k}]\ell(z)^k)\bigr)\ge-(m-j),$$
or equivalently
$$\nu\bigl(\sum_k\tfrac{2^k}{(k-j)!}([z^{m-k}]\ell(2z)^k)\bigr)\ge j\text{ with equality iff }\tbinom j{m-j}\text{ is odd.}$$
Since $\ell(2z)\equiv 1+z$ mod 2, and $\nu((k-j)!)\le k-j$ with equality iff $k=j$, all terms in the sum have $\nu(-)\ge j$ with equality iff
$k=j$ and $\binom j{m-j}$ is odd.\end{proof}

\begin{lem}\label{qblem} Let $q_m(-)$ be as in Definition \ref{qdef}, and let $x$ be any integer. Then
$\nu(q_m(x))\ge \nu(m!)$ with equality iff $\binom {x-m}m$ is odd.\end{lem}
\begin{proof} We have
$$q_m(x)=(x+1)_m-\sum_{j=1}^{m-1}2^{m-j}b_{j,m}q_j(x).$$
Note that $\nu((x+1)_m)\ge\nu(m!)$ with equality iff $\binom{x+1}m$ is odd. By induction, the $j$-term $T_j$ in the sum satisfies
$$\nu(T_j)\ge m-j+\nu(m!)-\nu(j!)-(m-j)+\nu(j!)=\nu(m!)$$
with equality iff $\binom j{m-j}$ is odd and $\binom{x-j}j$ is odd.
This implies the inequality. Equality occurs iff
\begin{equation}\label{bqsum}\tbinom{x+1}m+\sum_{j=0}^{m-1}\tbinom j{m-j}\tbinom{x-j}j\end{equation}
is odd.
By Lemma \ref{mod4lem}, $\displaystyle{\sum_{j=0}^m\tbinom j{m-j}\tbinom{x-j}j\equiv\tbinom{x+1}m}$ mod 2. Thus the expression in (\ref{bqsum})
is congruent to $\binom{x-m}m$, establishing the claim.
\end{proof}

Proposition \ref{refine} follows immediately from Theorem \ref{qthm2} and the following result, which is a refinement of Lemma \ref{qblem}.
\begin{thm} \label{refthm}If $m$ is a positive integer and $x$ is any integer, then, mod $4$,
$$q_m(x)/m!\equiv\tbinom{x-m}m+\begin{cases}2\tbinom{x-m}{m-2}&\text{if $x$ and $m$ are even}\\
0&\text{otherwise.}\end{cases}$$
\end{thm}

The proof of Theorem \ref{refthm} requires several subsidiary results.
\begin{lem} \label{mod4lem} If $m$ and $x$ are integers with $m\ge0$, then
$$\sum_{j=0}^m\tbinom j{m-j}\tbinom{x-j}j\equiv\tbinom{x+1}m+2\tbinom{x+1}{m-1}\pmod4.$$
\end{lem}
\begin{proof} This follows easily from Jensen's Formula (see e.g., \cite{G}), which says that if $A$, $B$, and $D$ are integers with $D\ge0$, then
$$\sum_{j=0}^D\tbinom{j+B}{D-j}\tbinom{A-j}j=\sum_{j=0}^D(-1)^j\tbinom{A+B-j}{D-j}.$$
This implies that the sum in our lemma equals $\displaystyle{\sum_{j=0}^m(-1)^j\tbinom{x-j}{m-j}}$. We prove that this is congruent, mod 4,
to the RHS of our lemma when $x\ge0$ by induction on $x$. The formula is easily seen to be true if $x=0$ (note that when $x=0$ and $m=1$ the LHS equals $-1$ while the RHS equals 3), and the induction step is by Pascal's formula. For $x<0$, let $y=-x$ with $y>0$. The equation to be proved becomes
$$\sum_{j=0}^m\tbinom{m+y-1}{m-j}\equiv\tbinom{y+m-2}m-2\tbinom{y+m-3}{m-1}\pmod4.$$
When $y=1$, both sides equal $\delta_{m,0}+2\delta_{m,1}$ and the result follows by induction on $y$ using Pascal's formula.
\end{proof}

The next result refines Lemma \ref{blem}.
\begin{lem}\label{bjmlem} If $b_{j,m}$ is as in Theorem \ref{qdef}, then, mod 4,
$$2^{m-j}j!b_{j,m}/m!\equiv\tbinom j{m-j}+2c_{j,m},\text{ where }c_{j,m}=\begin{cases}\tbinom j{m-j-1}&\text{if $j$ is even}\\
\tbinom j{m-j-2}&\text{if $j$ is odd.}\end{cases}$$
\end{lem}
\begin{proof} As in the proof of \ref{blem}, we have
\begin{equation}\label{b4eq}2^{m-j}j!b_{j,m}/m!=\sum_{k\ge j}\tfrac{2^{k-j}}{(k-j)!}([z^{m-k}]\ell(2z)^k).\end{equation}
Since, mod 4, $\ell(2z)\equiv 1-z-2z^3$, and $2^{k-j}/(k-j)!\equiv0$ unless $k-j$ equals 0 or a 2-power,  (\ref{b4eq})
equals \begin{eqnarray*}&&[z^{m-j}](1-z-2z^3)^j+2\sum_{e\ge0}[z^{m-j-2^e}](1-z-2z^3)^{j+2^e}\\
&\equiv&[z^{m-j}](1-z-2z^3)^j+2\sum_{e\ge0}\binom{j+2^e}{m-j-2^e}.\end{eqnarray*}
Replace $m-j$ by $\ell$. We must prove, mod 4,
\begin{equation}A_{j,\ell}+2B_{j,\ell}\equiv C_{j,\ell}+2D_{j,\ell},\label{ABCD}\end{equation}
where
$$A_{j,\ell}=\tbinom j\ell,\quad C_{j,\ell}=[z^\ell](1-z-2z^3)^j,\quad D_{j,\ell}=\sum_{e\ge0}\tbinom{j+2^e}{\ell-2^e},$$
and
$$ B_{j,\ell}=\begin{cases}\tbinom j{\ell-1}&\text{$j$ even}\\ \tbinom j{\ell-2}&\text{$j$ odd.}\end{cases}$$
If $j=0$, both sides of (\ref{ABCD}) are congruent to $\delta_{\ell,0}+2\delta_{\ell,1}$. For the RHS, note that
if $\ell=2^f$ with $f\ge1$, then $2 D_{0,\ell}\equiv0$ as it obtains a 2 from $e=f$ and from $e=f-1$.

Having proved the validity of (\ref{ABCD}) when $j=0$, we proceed by induction on $j$. If $j$ is even, then, mod 4,
\begin{eqnarray*}&&A_{j+1,\ell}+2B_{j+1,\ell}-C_{j+1,\ell}-2D_{j+1,\ell}\\
&=&A_{j,\ell}+A_{j,\ell-1}+2(B_{j,\ell-1}+B_{j,\ell-2})-(C_{j,\ell}-C_{j,\ell-1}-2C_{j,\ell-3})\\
&&-2(D_{j,\ell}+D_{j,\ell-1})\\
&\equiv&(A_{j,\ell}+2B_{j,\ell}-C_{j,\ell}-2D_{j,\ell})+(A_{j,\ell-1}+2B_{j,\ell-1}-C_{j,\ell-1}-2D_{j-1,\ell})\\
&&-2B_{j,\ell}+2B_{j,\ell-2}+2C_{j,\ell-1}+2C_{j,\ell-3}\\
&\equiv&-2\tbinom j{\ell-1}+2\tbinom j{\ell-3}+2\tbinom j{\ell-1}+2\tbinom j{\ell-3}\equiv0,
\end{eqnarray*}
and a similar argument works when $j$ is odd.
\end{proof}

The following result relates the even parts in \ref{refthm} and \ref{bjmlem}.
\begin{lem} \label{pclem}Let
$$p_j(x)=\begin{cases}\tbinom{x-j}{j-2}&\text{$x$ and $j$ even}\\0&\text{otherwise}\end{cases}\text{ and }c_{j,m}=\begin{cases}\tbinom j{m-j-1}&
\text{$j$ even}\\ \tbinom j{m-j-2}&\text{$j$ odd.}\end{cases}$$
Then, mod 2, if $x$ and $m$ are integers with $m\ge0$,
\begin{equation}\label{pceq}\tbinom{x+1}{m-1}\equiv\sum_{j=1}^m\bigl(\tbinom j{m-j}p_j(x)+c_{j,m}\tbinom{x-j}j\bigr).\end{equation}
\end{lem}
\begin{proof} First let $x$ be odd. By Lemma \ref{mod4lem}, mod 2,
$$\tbinom{x+1}{m-1}\equiv \sum_{j}\tbinom j{m-j-1}\tbinom{x-j}j.$$
Since $p_j(x)=0$ and $\binom{x-j}j\equiv0$ for odd $j$, this is equivalent to (\ref{pceq}) in this case.

Now suppose $x$ is even and $m$ odd. We must prove, mod 2,
$$\tbinom{x+1}{m-1}\equiv\sum_{j\text{ odd}}\tbinom j{m-j-2}\tbinom{x-j}j+\sum_{j\text{ even}}\bigl(\tbinom j{m-j}\tbinom{x-j}{j-2}+\tbinom j{m-j-1}\tbinom{x-j}j\bigr).$$
By \ref{mod4lem}, the LHS is congruent to $\sum\binom j{m-j-1}\binom{x-j}j$. If $j$ is odd, $\binom j{m-j-1}\equiv\binom j{m-j-2}$, and if
$j$ is even, $\binom j{m-j}\equiv0$. The desired result is now immediate.

Finally suppose $x$ and $m$ are both even. Again using \ref{mod4lem}, we must show
$$\sum_{j\text{ odd}}\tbinom j{m-j-1}\tbinom{x-j}j\equiv \sum_{j\text{ even}}\tbinom j{m-j}\tbinom{x-j}{j-2}+\sum_{j\text{ odd}}\tbinom j{m-j-2}\tbinom{x-j}j$$
since $\binom j{m-j-1}\equiv0$ if $j$ is even.
The terms on the LHS  combine with the $j$-odd terms on the RHS to yield $\displaystyle{\sum_{j\text{ odd}}\tbinom{j+1}{m-j-1}\tbinom{x-j}j}$. Letting $k=j+1$, this becomes $\displaystyle\sum_{k\text{ even}}\tbinom k{m-k}\tbinom{x-k+1}{k-1}$.
Since $x$ and $k$ are even, $\binom{x-k+1}{k-1}\equiv\binom{x-k}{k-2}$, and so all terms cancel.
\end{proof}

Now we easily prove Theorem \ref{refthm}.
\begin{proof}[Proof of Theorem \ref{refthm}.] The proof is by induction on $m$, with the case $m=1$ immediate. Using notation of \ref{pclem},
equation (\ref{q2d}) yields, mod 4,
\begin{eqnarray*}q_m(x)/m!&=&\tbinom{x+1}m-\sum_{j=1}^{m-1}\frac{j!2^{m-j}b_{j,m}}{m!}\frac{q_j(x)}{j!}\\
&\equiv&-2\tbinom{x+1}{m-1}+\sum_{j=0}^m\tbinom j{m-j}\tbinom{x-j}j-\sum_{j=1}^{m-1}\bigl(\tbinom j{m-j}+2c_{j,m}\bigr)\bigl(\tbinom{x-j}j+2p_j(x)\bigr)\\
&\equiv&\tbinom{x-m}m-2\bigl(\tbinom{x+1}{m-1}-\sum_{j=1}^{m-1}\bigl(\tbinom j{m-j}p_j(x)+c_{j,m}\tbinom{x-j}j\bigr)\\
&\equiv&\tbinom{x-m}m+2p_m(x),\end{eqnarray*}
as desired. Here we have used \ref{mod4lem} and \ref{bjmlem} at the second step and \ref{pclem} at the last step.
\end{proof}

\section{Relationship with Hensel's Lemma}\label{Clsec}

In \cite{D2}, the author introduced Lemma \ref{D2lem} and applied it to study $\nu(T_5(-))$ and $\nu(T_6(-))$ similarly to what we do here for all
$T_n(-)$. Clarke was quick to observe in \cite{Cl} that if $T_n(-)$ is considered as a function on $\zt$, then our conclusion that
$\nu(T_n(x))=\nu(x-x_0)+c_0$ when $x$ is restricted to a congruence class $C$ can be interpreted as saying that $T_n(x_0)=0$. He showed that
if $T_n(x_0)=0$ and  $|T_n'(x_0)|\ne0$, then
$$|T_n(x)|=|x-x_0||T_n'(x_0)|$$
on a neighborhood of $x_0$, which corresponds to our congruence class $C$. Here again $|x|=1/2^{\nu(x)}$ on $\zt$, and $d(x,y)=|x-y|$ defines the metric. Also, $T_n'$ denotes the derivative.
Moreover, Clarke noted that the
iteration toward the root $x_0$ in our theorems is a disguised form of Hensel's Lemma for convergence toward a root of the function $T_n$.

We illustrate by considering the root  of $T_{13}$ of the form $4x_0+1$. See Theorem \ref{36} and Table \ref{bigtbl}.
For our iteration toward $x_0$, let \begin{equation}\label{g2}g(x)=\nu(T_{13}(4x+1))-10.\end{equation}
Then $g(0)=1$, $g(2^1)=5$, $g(2^1+2^5)=6$, etc. Thus our early approximation to $4x_0+1$
is \begin{equation}\label{156}1+4(2^1+2^5+2^6),\end{equation}
and, continuing, we obtain that the last 18 digits in the binary expansion of $4x_0+1$ are \begin{equation}\label{binary}111001001110001001.\end{equation}
Note that each 1 in the binary expansion requires a separate calculation.

Now we describe the Hensel point of view, following Clarke (\cite{Cl}).
He showed that $$T_n'(k)=\sum\tbinom n{2i+1}(2i+1)^kL(2i+1),$$ where $$L(2i+1)=\dstyle{\sum_{j=1}^\infty(-1)^{j-1}(2i)^j/j}$$
is the 2-adic logarithm.
Hensel's Lemma applied to an analytic function $f$ involves the iteration $k_{n+1}=k_n-\frac{f(k_n)}{f'(k_n)}$, which, under favorable
hypotheses, converges to a root of $f$. We have $f=T_{13}$. Using {\tt Maple}, we find
\begin{equation}\label{T13p}\nu(T_{13}'(k))=\begin{cases}8&k\equiv1,2\ (4)\\
9&k\equiv0\ (4)\\
\ge11&k\equiv3\ (4).\end{cases}\end{equation}
To {\bf prove} this, which involves an infinite sum (for $L$) and infinitely many values of $k$, first note that our only claim is
about the mod $2^{11}$ value of $T_{13}'$, and so the sum for $L$ may be stopped after $j=12$. Since $L(2i+1)\equiv0$ mod 4, we are only concerned
with $(2i+1)^k$ mod $2^9$. Since $(2i+1)^k$ mod $2^9$ has period $2^8$ in $k$, performing the computation for 256 values of $k$ would suffice.

Let $k_0=1$. Then {\tt Maple} computes that $k_1=1-\tfrac{T_{13}(1)}{T'_{13}(1)}$ has binary expansion ending 1001001, and so agrees with (\ref{binary}) mod 64. Next $k_2=k_1-\frac{T_{13}(k_1)}{T_{13}'(k_1)}$ has binary expansion ending 0001110001001, agreeing with (\ref{binary}) mod $2^{12}$. Finally $k_3=k_2-\frac{T_{13}(k_2)}{T_{13}'(k_2)}$ agrees with (\ref{binary}), and hence is correct at least mod $2^{18}$. This is much faster convergence than ours.

Let $\theta(x)=T_{13}(4x+1)$. Our algorithm essentially applies Hensel's Lemma to $\theta(x)$, but just takes the leading term each time.
For all $x$, $\nu(\theta'(x))=\nu(4T_{13}'(4x+1))=10$, and so our $g(x)$ equals $\nu(\theta(x)/\theta'(x))$. Thus when we let $x_{i+1}=x_i+2^{g(x_i)}$,
we are adding the leading term of $\theta(x_i)/\theta'(x_i)$. Once the limiting value, which we denote by $x_0$, is found, the root of $T_{13}$ is $4x_0+1$.

In \cite{Cl}, Clarke defines, for an analytic function $f$,
$$\bg(x,h)=\frac{f(x+h)-f(x)-hf'(x)}{h^2}$$
and shows that if  $f(x_0)=0$ and $|\bg(x,h)|\le2^r$ for all relevant $x$ and $h$, then the desired formula
$$|f(x)|=|x-x_0||f'(x_0)|$$ holds for all $x$ satisfying
\begin{equation}\label{gin}|x-x_0|<|f'(x_0)|/2^r.\end{equation}
He also notes that our $T_n(-)$ are analytic when restricted to all 2-adic integers of either parity.

For $f=T_{13}$, {\tt Maple} suggests that $\nu(\bg(x,h))\ge7$ for $h$ even, with equality iff $x+h\equiv0,3$ mod 4. This can be easily
proved using {\tt Maple} calculations and some elementary arguments. Hence $r=-7$. Using (\ref{T13p}), we obtain the results of Table \ref{bigtbl} for $n=13$
which are
listed there as $C=[1,2\ (4)]$ and $[0,4\ (8)]$, since
$$\frac{|T_{13}'(x_0)|}{2^{-7}}=\begin{cases}2^{-1}&x_0\equiv1,2\ (4)\\
2^{-2}&x_0\equiv0\ (4).\end{cases}$$
Being less than this requires $|x-x_0|\le2^{-2}$ or $2^{-3}$ in (\ref{gin}), whose reciprocals are the moduli of the congruence classes in
Table \ref{bigtbl}.

Another result of \cite{Cl} gives a condition,
\begin{equation}\label{convcond}|f(x)|<\min\bigl(\tfrac{|f'(x)|}{2^{k-1}},\tfrac{|f'(x)|^2}{2^r}\bigr)\end{equation}
(where $|\bg(x,h)|\le2^r$ and $f$ is analytic on $c+2^k\zt$),
which guarantees that iteration of Hensel from $x$ converges to a root of $f$. For
$f=T_{13}$ and $x\equiv3\ (4)$,  $|T_{13}'(x)|^2/2^r\le(2^{-11})^2/2^{-7}=2^{-15}$ by (\ref{T13p}), while by Table \ref{bigtbl} $|T_{13}(x)|$ takes on values $2^{-11}$, $2^{-13}$, and $2^{-15}$. Thus the condition (\ref{convcond}) does not
hold, consistent with our finding in Table \ref{bigtbl} that $|T_{13}(x)|$ is constant on balls about 7, 3, and 11, so there is no root in these
neighborhoods.

Clarke's approach is a very attractive alternative to ours. It converges faster, and it is more closely associated
with analytic methods, such as the Hensel/Newton convergence algorithm. On the other hand, there is a certain combinatorial
simplicity to our approach, especially Lemma \ref{D2lem} and its reduction to consideration of expressions such as (\ref{term}) and
(\ref{t3e2}), and subsequently to (\ref{alpheq}).
 We  find it very attractive that for each $f=T_n$, it seems likely
that $\zt$ can be partitioned into finitely many balls $B(x_0,\eps)$ on each of which $|f(x)|$ is linear in $|x-x_0|$ (including the possibility
that it is constant). It is not clear which approach will be the better way to establish this.

\def\line{\rule{.6in}{.6pt}}

\end{document}